\input amstex
\documentstyle{amsppt}
\magnification=\magstep1                        
\hsize6.5truein\vsize8.9truein                  
\NoRunningHeads
\loadeusm

\magnification=\magstep1                        
\hsize6.5truein\vsize8.9truein                  
\NoRunningHeads
\loadeusm

\document
\topmatter

\title
The asymptotic distance between an ultraflat unimodular polynomial and its conjugate reciprocal
\endtitle

\rightheadtext{ultraflat unimodular unimodular polynomials}

\author Tam\'as Erd\'elyi
\endauthor

\address Department of Mathematics, Texas A\&M University,
College Station, Texas 77843, College Station, Texas 77843 \endaddress

\endthanks


\date February 10, 2019
\enddate

\email terdelyi\@math.tamu.edu
\endemail

\dedicatory Dedicated to the memory of Jean-Pierre Kahane \enddedicatory

\abstract
Let
$${\Cal K}_n := \left\{p_n: p_n(z) = \sum_{k=0}^n{a_k z^k}, \enskip  a_k \in {\Bbb C}\,,
\enskip |a_k| = 1 \right\}\,.$$
A sequence $(P_n)$ of polynomials $P_n \in {\Cal K}_n$ is called ultraflat if 
$(n + 1)^{-1/2}|P_n(e^{it})|$ converge to $1$ uniformly in $t \in {\Bbb R}$.
In this paper we prove that  
$$\frac{1}{2\pi} \int_0^{2\pi}{\left| (P_n - P_n^*)(e^{it}) \right|^q \, dt} \sim 
\frac{{2}^q \Gamma \left(\frac{q+1}{2} \right)}{\Gamma \left(\frac q2 + 1 \right) \sqrt{\pi}} \,\, n^{q/2}$$ 
for every ultraflat sequence $(P_n)$ of polynomials $P_n \in {\Cal K}_n$ and for every $q \in (0,\infty)$, where
$P_n^*$ is the conjugate reciprocal polynomial associated with $P_n$, $\Gamma$ is the usual gamma function, and
the $\sim$ symbol means that the ratio of the left and right hand sides converges to $1$ as $n \rightarrow \infty$.
Another highlight of the paper states that 
$$\frac{1}{2\pi}\int_0^{2\pi}{\left| (P_n^\prime - P_n^{*\prime})(e^{it}) \right|^2 \, dt} \sim \frac{2n^3}{3}$$
for every ultraflat sequence $(P_n)$ of polynomials $P_n \in {\Cal K}_n$.
We prove a few other new results and reprove some interesting old results as well. 
\endabstract
\endtopmatter

\head 1. Introduction \endhead
Let
$${\Cal K}_n := \left\{p_n: p_n(z) = \sum_{k=0}^n{a_k z^k}, \enskip  a_k \in {\Bbb C}\,,
\enskip |a_k| = 1 \right\}\,.$$
The class ${\Cal K}_n$ is often called the collection of all (complex) unimodular polynomials of degree $n$.
Let
$${\Cal L}_n := \left\{p_n: p_n(z) = \sum_{k=0}^n{a_k z^k}, \enskip  a_k \in \{-1,1\} \right\}\,.$$
The class ${\Cal L}_n$ is often called the collection of all (real) unimodular polynomials of degree $n$.
By Parseval's formula,
$$\int_{0}^{2\pi}{\left| P_n(e^{it}) \right|^2 \, dt} = 2\pi(n+1)$$
for all $P_n \in {\Cal K}_n$. Therefore
$$\min_{t \in {\Bbb R}}{|P_n(e^{it})|} \leq \sqrt{n + 1} \leq \max_{t \in \Bbb R}{|P_n(e^{it})|}\,.$$
An old problem (or rather an old theme) is the following. 

\proclaim{Problem 1.1 (Littlewood's Flatness Problem)}
How close can a polynomial $P_n \in {\Cal K}_n$ or $P_n \in {\Cal L}_n$ come to satisfying
$$|P_n(e^{it})| = \sqrt{n + 1}\,, \qquad t \in {\Bbb R}? \tag 1.1$$
\endproclaim

Obviously (1.1) is impossible if $n \geq 1$. So one must
look for less than (1.1), but then there are various ways of seeking
such an ``approximate situation". One way is the following.
In his paper [Li1] Littlewood had suggested that, conceivably, there might exist a sequence
$(P_n)$ of polynomials $P_n \in {\Cal K}_n$ (possibly even $P_n \in {\Cal L}_n$) such
that $(n + 1)^{-1/2}|P_n(e^{it})|$ converge to $1$ uniformly in $t \in {\Bbb R}$.
We shall call such sequences of unimodular polynomials ``ultraflat". More precisely, we give the
following definition.

\proclaim{Definition 1.2} Given a positive number $\varepsilon$,
we say that a polynomial $P_n \in {\Cal K}_n$ is $\varepsilon$-flat if
$$(1 - \varepsilon)\sqrt{n + 1} \leq |P_n(e^{it})| \leq (1 + \varepsilon)\sqrt{n + 1}\,,
\qquad t \in {\Bbb R}\,.$$
\endproclaim

\proclaim{Definition 1.3} Let $(n_k)$ be an increasing sequence of positive integers. Given a sequence 
$(\varepsilon_{n_k})$ of positive numbers tending to $0$, we say that a sequence $(P_{n_k})$ of polynomials 
$P_{n_k} \in {\Cal K}_{n_k}$ is $(\varepsilon_{n_k})$-ultraflat if each $P_{n_k}$ is $(\varepsilon_{n_k})$-flat. 
We simply say that a sequence $(P_{n_k})$ of polynomials $P_{n_k} \in {\Cal K}_{n_k}$ is
ultraflat if it is $(\varepsilon_{n_k})$-ultraflat with a suitable sequence $(\varepsilon_{n_k})$ 
of positive numbers tending to $0$.
\endproclaim

The existence of an ultraflat sequence of unimodular polynomials seemed very unlikely, in view of a
1957 conjecture of P. Erd\H os (Problem 22 in [Er]) asserting that, for all
$P_n \in {\Cal K}_n$ with $n \geq 1$,
$$\max_{t \in {\Bbb R}}{|P_n(e^{it})|} \geq (1 + \varepsilon) \sqrt{n+1}\,, \tag 1.2$$
where $\varepsilon > 0$ is an absolute constant (independent of $n$).
Yet, refining a method of K\"orner [K\"o], Kahane [Ka] proved that there exists
a sequence $(P_n)$ with $P_n \in {\Cal K}_n$ which is $(\varepsilon_n)$-ultraflat, where
$\varepsilon_n = O\left(n^{-1/17} \sqrt{\log n} \right)\,.$
(Kahane's paper contained though a slight error which was corrected in [QS2].)  
Thus the Erd\H os conjecture (1.2) was disproved for the classes ${\Cal K}_n$.
For the more restricted class ${\Cal L}_n$ the analogous Erd\H os conjecture
is unsettled to this date. It is a common belief that the analogous Erd\H os conjecture
for ${\Cal L}_n$ is true, and consequently there is no ultraflat sequence of polynomials
$P_n \in {\Cal L}_n$.
An interesting result related to Kahane's breakthrough is given in [Be].
For an account of some of the work done till the mid 1960's, see Littlewood's book [Li2]
and [QS2]. 

If $Q_n$ is a polynomial of degree $n$ of the form
$$Q_n(z) = \sum_{k=0}^n {a_k z^k}\,, \qquad a_k \in {\Bbb C}\,,$$
then its conjugate polynomial is defined by
$$Q_n^*(z) := z^n \overline{Q}_n(1/z) := \sum_{k=0}^n {\overline{a}_{n-k}z^k}\,.$$

Let $(\varepsilon_n)$ be a sequence of positive numbers tending to $0$.
Let the sequence $(P_n)$ of polynomials $P_n \in {\Cal K}_n$ be $(\varepsilon_{n})$-ultraflat. 
We write
$$P_n(e^{it}) = R_n(t) e^{i\alpha_n(t)}\,, \qquad R_n(t) = |P_n(e^{it})|\,, \qquad t \in {\Bbb R} \,. \tag 1.3$$
It is simple to show that $\alpha_n$ can be chosen to be in $C^\infty({\Bbb R})$.
This is going to be our understanding throughout the paper.
It is easy to find a formula for $\alpha_n(t)$ in terms of $P_n$. We have 
$$\alpha_n^{\prime}(t) = \text {\rm Re} \left( \frac{e^{it}P_n^{\prime}(e^{it})}{P_n(e^{it})} \right) \,, \tag 1.4$$ 
see formulas (7.1) and (7.2) on p.~564 and (8.2) on p.~565 in [Sa1].
The angular function $\alpha_n^*$ and modulus function $R_n^*=R_n$ associated with
the polynomial $P_n^*$ are defined by
$$P_n^*(e^{it}) = R_n^*(t) e^{i\alpha_n^*(t)}\,, \qquad R_n^*(t) = |P_n^*(e^{it})|\,.$$ 
Similarly to $\alpha_n$, the angular function $\alpha_n^*$ can also be chosen to be in $C^\infty({\Bbb R})$ on ${\Bbb R}$. 
By applying formula (1.4) to $P_n^*$, it is easy to see that
$$\alpha_n^\prime(t) + {\alpha_n^*}^\prime(t) = n\,, \qquad t \in {\Bbb R} \,. \tag 1.5$$

The structure of ultraflat sequences of unimodular polynomials is studied in
[Er1], [Er2], [Er3], and [Er4], where several conjectures of Saffari are proved. 
In [Er6], based on the results in [Er1], we proved yet another conjecture of  
Queffelec and Saffari, see (1.30) in [QS2]. Namely we proved asymptotic formulas 
for the $L_q$ norms of the real part and the derivative of the real part of ultraflat 
unimodular polynomials on the unit circle. A recent paper of Bombieri and Bourgain 
[BB] is devoted to the construction of ultraflat sequences of unimodular polynomials. 
In particular, they obtained a much improved estimate for the error term. A major 
part of their paper deals also with the long-standing problem of the effective 
construction of ultraflat sequences of unimodular polynomials.  

For $\lambda \ge 0$, let
$${\Cal K}_n^{\lambda} := \left\{P_n: P_n(z) = \sum_{k=0}^n{a_k k^\lambda z^k}, 
\enskip  a_k \in {\Bbb C}\,, \enskip |a_k| = 1 \right\} \,.$$
Ultraflat sequences $(P_n)$ of polynomials $P_n \in {\Cal K}_n^{\lambda}$ are defined and 
studied thoroughly in [EN] where various extensions of Saffari's conjectures have been proved. 
  
In [Er2] we examined how far an ultraflat unimodular polynomial is from being conjugate reciprocal, 
and we proved the following three theorems.  

\proclaim{Theorem 1.4} Let $(P_n)$ be an ultraflat sequence of polynomials $P_n \in {\Cal K}_n$. 
We have
$$\frac{1}{2\pi}\int_0^{2\pi}{\left( |P_n^\prime(e^{it})| - |P_n^{*\prime}(e^{it})| \right)^2 \, dt}
= \left( \frac 13 + \gamma_n \right) n^3\,,$$
where $(\gamma_n)$ is a sequence of real numbers converging to $0$.
\endproclaim

\proclaim{Theorem 1.5} Let $(P_n)$ be an ultraflat sequence of polynomials $P_n \in {\Cal K}_n$.
If the coefficients of $P_n$ are denoted by $a_{k,n}$, that is, 
$$P_n(z) = \sum_{k=0}^n{a_{k,n} z^k}\,, \qquad k=0,1,\ldots, n, \quad n=1,2,\ldots\,,$$
then
$$\sum_{k=0}^n{k^2 \left| a_{k,n} - \overline{a}_{n-k,n} \right|^2} = 
\frac{1}{2\pi}\int_0^{2\pi}{\left|(P_n^\prime - P_n^{*\prime})(e^{it})\right|^2 \, dt} 
\geq \left( \frac 13 + h_n \right) n^3\,,$$
where $(h_n)$ is a sequence of real numbers converging to $0$.
\endproclaim

\proclaim{Theorem 1.6} Let $(P_n)$ be an ultraflat sequence of polynomials $P_n \in {\Cal K}_n$.
Using the notation of Theorem 1.5 we have
$$\sum_{k=0}^n{\left| a_{k,n} - \overline{a}_{n-k,n} \right|^2} = 
\frac{1}{2\pi}\int_0^{2\pi}{\left| (P_n - P_n^*)(e^{it}) \right|^2 \, dt}
\geq \left( \frac 13 + h_n \right) n\,,$$
where $(h_n)$ is the same sequence of real numbers converging to $0$ as in Theorem 1.5.
\endproclaim

There are quite a few recent publications on or related to ultraflat sequences of unimodular 
polynomials. Some of them (not mentioned before) are are [Bo], [Sa2], [QS1], [Od], and [Mo]. 

\head 2. Results \endhead

Theorems 2.1--2.4 and 2.6 are new, Theorems 2.5 and 2.7 recapture old results.  

In our results below $\Gamma$ denotes the usual gamma function, and the $\sim$ symbol means
that the ratio of the left and right hand sides converges to $1$ as $n \rightarrow \infty$.

\proclaim{Theorem 2.1}
If $(P_n)$ is an ultraflat sequence of polynomials $P_n \in {\Cal K}_n$ and $q \in (0,\infty)$, then
$$\frac{1}{2\pi} \int_0^{2\pi}{\left| (P_n - P_n^*)(e^{it}) \right|^q \, dt} \sim 
\frac{{2}^q \Gamma \left(\frac{q+1}{2} \right)}{\Gamma \left(\frac q2 + 1 \right) \sqrt{\pi}} \,\, n^{q/2}\,.$$
\endproclaim 

Our next theorem is a special case ($q=2$) of Theorem 2.1.  Compare it with Theorem 1.6.

\proclaim{Theorem 2.2}
Let $(P_n)$ be an ultraflat sequence of polynomials $P_n \in {\Cal K}_n$.
If the coefficients of $P_n$ are denoted by $a_{k,n}$, that is,  
$$P_n(z) = \sum_{k=0}^n{a_{k,n} z^k}\,, \qquad k=0,1,\ldots, n, \quad n=1,2,\ldots\,,$$
then
$$\sum_{k=0}^n{\left| a_{k,n} - \overline{a}_{n-k,n} \right|^2} = 
\frac{1}{2\pi} \int_0^{2\pi}{\left| (P_n - P_n^*)(e^{it}) \right|^2 \, dt} \sim 2n\,.$$
\endproclaim

Our next theorem should be compared with Theorem 1.5.

\proclaim{Theorem 2.3} Let $(P_n)$ be an ultraflat sequence of polynomials $P_n \in {\Cal K}_n$.
Using the notation in Theorem 2.2 we have
$$\sum_{k=0}^n{k^2 \left| a_{k,n} - \overline{a}_{n-k,n} \right|^2} = 
\frac{1}{2\pi}\int_0^{2\pi}{\left| (P_n^\prime - P_n^{*\prime})(e^{it}) \right|^2 \, dt} \sim \frac{2n^3}{3}\,.$$
\endproclaim

We also prove the following result.

\proclaim{Theorem 2.4} 
If $(P_n)$ is an ultraflat sequence of polynomials $P_n \in {\Cal K}_n$ and $q \in (0,\infty)$, then
$$\frac{1}{2\pi} \int_0^{2\pi}{\left| \frac{d}{dt} |(P_n - P_n^*)(e^{it})| \right|^q \, dt} 
\sim \frac{\Gamma \left(\frac{q+1}{2} \right)}{(q+1)\Gamma \left(\frac q2 + 1 \right) \sqrt{\pi}} \,\, n^{3q/2}\,.$$
\endproclaim 

As a Corollary of Theorem 2.2 we can recapture Saffari's ``near orthogonality conjecture" raised in [Sa]
and proved first in [Er4].

\proclaim{Theorem 2.5}
Let $(P_n)$ be an ultraflat sequence of polynomials $P_n \in {\Cal K}_n$.
Using the notation in Theorem 2.2 we have
$$\sum_{k=0}^n{a_{k,n}a_{n-k,n}} = o(n)\,.$$
\endproclaim

As a Corollary of Theorem 2.3 we can easily prove a new ``near orthogonality" formula.

\proclaim{Theorem 2.6}
Let $(P_n)$ be an ultraflat sequence of polynomials $P_n \in {\Cal K}_n$.
Using the notation in Theorem 2.2 we have
$$\sum_{k=0}^n{k^2a_{k,n}a_{n-k,n}} = o(n^3)\,.$$
\endproclaim

Finally we recapture the asymptotic formulas for the real part and the derivative of the real part of 
ultraflat unimodular polynomials proved in [Er5] first.  

\proclaim{Theorem 2.7}
If $(P_n)$ is an ultraflat sequence of unimodular polynomials $P_n \in {\Cal K}_n$,
and $q \in (0,\infty)$, then for $f_n(t) := \text{\rm Re}(P_n(e^{it}))$ we have
$$\frac{1}{2\pi} \int_0^{2\pi}{\left| f_n(t) \right|^q \, dt} 
\sim \frac{\Gamma \left(\frac{q+1}{2} \right)}{\Gamma \left(\frac q2 + 1 \right) \sqrt{\pi}} \,\, n^{q/2}$$
and
$$\frac{1}{2\pi} \int_0^{2\pi}{\left| f_n^\prime(t) \right|^q \, dt}
\sim \frac{\Gamma \left(\frac{q+1}{2} \right)}{(q+1)\Gamma \left(\frac q2 + 1 \right) \sqrt{\pi}} \,\, n^{3q/2}\,.$$
\endproclaim

We remark that trivial modifications of the proof of Theorem 2.1--2.7 yield that 
the statement of the above theorem remains true if the ultraflat sequence $(P_n)$ of 
polynomials $P_n \in {\Cal K}_n$ is replaced by an ultraflat sequence
$(P_{n_k})$ of polynomials $P_{n_k} \in {\Cal K}_{n_k}$, $0 < n_1 < n_2 < \ldots$. 

As a Corollary of Theorem 2.2 we can recapture Saffari's ``near orthogonality conjecture" raised in [Sa] 
and proved first in [Er4].

\head 3. Lemmas \endhead 

To prove Theorems 2.1 and 2.2 we need a few lemmas. The first two are from [Er1]. 

\proclaim{Lemma 3.1 (Uniform Distribution Theorem for the Angular Speed)}
Suppose $(P_n)$ is an ultraflat sequence of polynomials $P_n \in {\Cal K}_n$. Then, with the notation (1.3), 
in the interval $[0,2\pi]$, the distribution of the normalized
angular speed $\alpha_n^\prime(t)/n$ converges to the uniform distribution
as $n \rightarrow \infty$. More precisely, we have
$$\text{\rm meas} (\{t \in [0,2\pi]: 0 \leq \alpha_n^\prime(t) \leq nx\}) = 2\pi x + \gamma_n(x)$$
for every $x \in [0,1]$, where 
$\displaystyle{\lim_{n \rightarrow \infty}{\max_{x \in [0,1]}{|\gamma_n(x)|}} = 0}$.   
\endproclaim

Our next lemma is a simple observation of Saffari [Sa1], which follows from (1.4), Bernstein's inequality, and 
the ultraflatness property given by Definition 1.3.

\proclaim{Lemma 3.2}
Suppose $(P_n)$ is an ultraflat sequence of polynomials $P_n \in {\Cal K}_n$. Then, with the notation (1.3),
we have
$$o_nn \leq \alpha_n^\prime(t) \leq n - o_nn\,, \qquad t \in {\Bbb R}\,, \tag 3.1$$
with real numbers $o_n$ converging to $0$.
\endproclaim

\proclaim{Lemma 3.3 (Negligibility Theorem for Higher Derivatives)}
Suppose $(P_n)$ is an ultraflat sequence of polynomials $P_n \in {\Cal K}_n$. Then, with the notation (1.3),
for every integer $r \geq 2$, we have
$$\max_{0 \leq t \leq 2\pi}{|\alpha_n^{(r)}(t)|} \leq \gamma_{n,r}n^r$$
with real numbers $\gamma_{n,r} > 0$ converging to $0$ for every fixed $r = 2,3,\ldots$.
\endproclaim

Our next lemma is a special case of Lemma 4.2 from [Er1].

\proclaim{Lemma 3.4}
Suppose $(P_n)$ is an ultraflat sequence of polynomials $P_n \in {\Cal K}_n$. 
Using notation (1.3), we have
$$\max_{0 \leq t \leq 2\pi}{|R_n^\prime(t)|} = \varphi_n n^{3/2}\,,
\qquad m = 1,2,\ldots\,,$$
with real numbers $\varphi_n$ converging to $0$. 
\endproclaim 

The next lemma follows simply from the ultraflatness property given by Definition 1.3.

\proclaim{Lemma 3.5} Let $q \in (0,\infty)$. Using the notation (1.3) and 
$$\beta_n(t) := \frac 12 (\alpha_n(t) - \alpha_n^*(t)) = \alpha_n(t) - \frac{nt}{2} - t_0$$ 
we have
$$\frac{1}{2\pi} \int_0^{2\pi}{\left| (P_n - P_n^*)(e^{it}) \right|^q \, dt} = 
\int_0^{2\pi}{\left| n^{1/2}(1 + \delta_n(t)) 2\sin(\beta_n(t)) \right|^q \, dt}$$
with real numbers $\delta_n(t)$ converging to $0$ uniformly on $[0,2\pi]$.  
\endproclaim

\proclaim{Lemma 3.5*}
Let $(P_n)$ be an ultraflat sequence of unimodular polynomials $P_n \in {\Cal K}_n$,
$q \in (0,\infty)$, and $f_n(t) := \text{\rm Re}(P_n(e^{it}))$. Using the notation (1.3) we have
$$\frac{1}{2\pi} \int_0^{2\pi}{\left| f_n(t) \right|^q \, dt} = 
\int_0^{2\pi}{\left| n^{1/2}(1 + \delta_n(t)) \cos(\alpha_n(t)) \right|^q \, dt}$$
with real numbers $\delta_n(t)$ converging to $0$ uniformly on $[0,2\pi]$.
\endproclaim

The next lemma follows simply from the ultraflatness property given by Definition 1.3 and Lemma 3.4.

\proclaim{Lemma 3.6} Let $q \in (0,\infty)$. Using the notation (1.3) and
$$\beta_n(t) := \frac 12 (\alpha_n(t) - \alpha_n^*(t)) = \alpha_n(t) - \frac{nt}{2} - t_0$$
we have
$$\frac{1}{2\pi} \int_0^{2\pi}{\left| \frac{d}{dt} |(P_n - P_n^*)(e^{it})| \right|^q \, dt}$$   
$$= \, \int_0^{2\pi}{\left| n^{1/2}(1 + \delta_n(t))2\cos(\beta_n(t)) \beta_n^\prime(t) + \eta_n(t) n^{3/2} \right|^q \,dt}$$
with real numbers $\delta_n(t)$ and $\eta_n(t)$ converging to $0$ uniformly on $[0,2\pi]$.
\endproclaim

\proclaim{Lemma 3.6*}
Let $(P_n)$ be an ultraflat sequence of unimodular polynomials $P_n \in {\Cal K}_n$,
$q \in (0,\infty)$, and $f_n(t) := \text{\rm Re}(P_n(e^{it}))$. Using the notation (1.3) we have
$$\frac{1}{2\pi} \int_0^{2\pi}{\left| f_n^\prime(t) \right|^q \, dt}$$ 
$$= \, \int_0^{2\pi}{\left| n^{1/2}(1 + \delta_n(t))\sin(\alpha_n(t)) \alpha_n^\prime(t) + \eta_n(t) n^{3/2} \right|^q \,dt}$$
with real numbers $\delta_n(t)$ converging to $0$ uniformly on $[0,2\pi]$.
\endproclaim

To prove Lemmas 3.8 and 3.9 we need the technical lemma below that follows by a simple calculation 
using formulas (6.2.1), (6.2.2), and (6.1.8) on pages 258 and 255 in [AS].

\proclaim{Lemma 3.7}
Assume that $A,B \in {\Bbb R}$, $B \neq 0$, $q > 0$, and $I \subset [0,2\pi]$ is an interval.
Then
$$\int_I{\left| \cos(Bt + A) \right|^q \,dt} = K(q) \text{\rm meas}(I) + \delta_1(A,B,q)$$ 
and
$$\int_I{\left| \sin(Bt + A) \right|^q \,dt} = K(q) \text{\rm meas}(I) + \delta_2(A,B,q)\,,$$ 
where
$$K(q) := \frac{1}{2\pi} \int_0^{2\pi}{\left| \sin t \right|^q \,dt} = 
\frac{\Gamma \left(\frac{q+1}{2} \right)}{\Gamma \left(\frac q2 + 1 \right) \sqrt{\pi}}$$
and 
$$\left| \delta_1(A,B,q) \right| \leq \pi \left| B \right|^{-1} \quad \text{and} \quad 
\left| \delta_2(A,B,q) \right| \leq \pi \left| B \right|^{-1}\,.$$ 
\endproclaim

Our final couple of lemmas take care of the most difficult part of the proof of Theorem 2.1.

\proclaim{Lemma 3.8}
Suppose that $\beta_n$, $n=1,2,\ldots$, are real-valued functions defined on $[0,2\pi]$ such that their second 
derivatives $\beta_n^{\prime\prime}$ are continuous on $[0,2\pi]$. Suppose also that  
$$\text{\rm meas} (\{t \in [0,2\pi]: |2\beta_n^\prime(t)| \leq nx\}) = \gamma(x) + \gamma_n(x)\,, 
\qquad x \in [0,1]\,, \tag 3.2$$  
where
$$\lim_{x \rightarrow 0+}{\gamma(x)} = \lim_{n \rightarrow \infty}{\max_{x \in [0,1]}{|\gamma_n(x)|}} = 0\,,\tag 3.3$$ 
and
$$\max_{0 \leq t \leq 2\pi}{|\beta_n^{\prime\prime}(t)|} \leq \gamma_{n,2}n^2 \tag 3.4$$
with real numbers $\gamma_{n,2} > 0$ converging to $0$. 
Then
$$\frac{1}{2\pi} \int_0^{2\pi}{\left| \sin(\beta_n(t)) \right|^q \,dt} \sim K(q) := 
\frac{\Gamma \left(\frac{q+1}{2} \right)}{\Gamma \left(\frac q2 + 1 \right) \sqrt{\pi}}\,. \tag 3.5$$
\endproclaim

\proclaim{Lemma 3.9}
Suppose that $\beta_n$, $n=1,2,\ldots$, are real-valued functions defined on $[0,2\pi]$ such that their second
derivatives $\beta_n^{\prime\prime}$ are continuous on $[0,2\pi]$ and 
$$\text{\rm meas} (\{t \in [0,2\pi]: |2\beta_n^\prime(t)| \leq nx\}) = 2\pi x + \gamma_n(x)\,, 
\qquad x \in [0,1]\,,  \tag 3.6$$ 
where (3.3) holds (with $\gamma(x) := 2\pi x$). Suppose also that (3.4) holds and  
$$\max_{t \in [0,2\pi]}{|\beta_n^\prime(t)|} \leq cn \tag 3.7$$
with an absolute constant $c > 0$. Then
$$\frac{1}{2\pi} \int_0^{2\pi}{\left| \cos(\beta_n(t)) n^{-1}\beta_n^\prime(t) \right|^q \,dt} 
\sim \frac{K(q)}{{2}^q(q+1)} := 
\frac{\Gamma \left(\frac{q+1}{2} \right)}{{2}^q(q+1)\Gamma \left(\frac q2 + 1 \right) \sqrt{\pi}}\,. \tag 3.8$$
\endproclaim

We note that conditions (3.6), (3.3), and (3.7) imply in a standard fashion that
$$\frac{1}{2\pi}\int_0^{2\pi}{\left| 2\beta_n^\prime(t) \right|^q\,dt} = \frac{n^{q}}{q+1} + \delta_{n,q}n^q \tag 3.9$$
with real numbers $\delta_{n,q}$ converging to $0$ for every fixed $q > 0$.

\demo{Proof of Lemma 3.8}
Let $\varepsilon > 0$ be fixed. Let $L_n := \gamma_{n,2}^{-1/4}\,.$  
We divide the interval $[0,2\pi]$ into subintervals 
$$I_m := [a_{m-1},a_m) := \left[\frac{(m-1)L_n}{n}, \frac{mL_n}{n}\right), \qquad m=1,2, \ldots, 
N - 1 := \left \lfloor \frac{2\pi n}{L_n} \right \rfloor\,,$$
and
$$I_N := [a_{N-1}, a_N) := \left[\frac{(N-1)L_n}{n}, 2\pi\right)\,.$$ 
For the sake of brevity let
$$A_{m-1} := \beta_n(a_{m-1})\,, \qquad m=1,2, \ldots, N\,,$$
and
$$B_{m-1} := \beta_n^\prime(a_{m-1})\,, \qquad m=1,2, \ldots, N\,.$$
Using Taylor's Theorem and assumption (3.4) we obtain that 
$$\left| \beta_n(t) - (A_{m-1} + B_{m-1}(t - a_{m-1}) \right| \leq \gamma_{n,2} n^2 (L_n/n)^2 
\leq \gamma_{n,2} \gamma_{n,2}^{-1/2} \leq \gamma_{n,2}^{1/2}$$
for every $t \in I_m$, where $\displaystyle{\lim_{n \rightarrow \infty}{\gamma_{n,2}^{1/2}} = 0}$. 
Hence the functions 
$$G_{n,q}(t) := \left\{\aligned 
|\sin&(A_{0} + B_{0}(t - a_{0}))|^q, \\ 
|\sin&(A_{1} + B_{1}(t - a_{0}))|^q, \\ 
& \vdots \\
|\sin&(A_{N-1} + B_{N-1}(t - a_{N-1}))|^q, \endaligned \right. \qquad
\aligned 
t & \in I_1\,, \\
t & \in I_2\,, \\ 
& \vdots \\ 
t & \in I_N\,, \\ \endaligned$$ 
and
$$F_{n,q}(t) := \left| \sin(\beta_n(t)) \right|^q$$
satisfy 
$$\lim_{n \rightarrow \infty}{\sup_{t \in [0,2\pi)}{| G_{n,q}(t) - F_{n,q}(t)|}} = 0\,. \tag 3.10$$
Therefore, in order to prove (3.5), it is sufficient to prove that
$$\int_0^{2\pi}{G_{n,q}(t) \,dt} \sim 2\pi K(q)\,. \tag 3.11$$
If $|B_{m-1}| \geq n\varepsilon$, then Lemma 3.7 gives
$$\left|\int_{I_m}{G_{n,q}(t) \,dt} - K(q)\text{\rm meas}(I_m)\right| \leq \frac{\pi}{n\varepsilon}\,.$$
By assumption (3.4) we have 
$\displaystyle{\lim_{n \rightarrow \infty}{L_n} = \lim_{n \rightarrow \infty}{\gamma_{n,2}^{-1/4}} = \infty}$, 
and hence  
$$\split \left| \sum_{m}{\int_{I_m}{G_{n,q}(t) \,dt}} - K(q)\sum_{m}{\text{\rm meas}(I_m)} \right| 
\leq N\frac{\pi}{n\varepsilon} 
& \leq \left(\frac{2\pi n}{L_n} + 1\right) \frac{\pi}{n\varepsilon} \cr 
& \leq \eta_n(\varepsilon)\,, \cr \endsplit \tag 3.12$$
where the summation is taken over all $m=1,2, \ldots, N$ for which $|B_{m-1}| \geq n\varepsilon$,
and where $(\eta_n(\varepsilon))$ is a sequence of real numbers tending to $0$.
Now let
$$E_{n,\varepsilon} := \bigcup_{m: \, |B_{m-1}| \leq n\varepsilon}{I_m}\,.$$
If $|B_{m-1}| \leq n\varepsilon$, then by assumption (3.4) we have
$$|\beta_n^\prime(t)| \leq |B_{m-1}| + \frac{L_n}{n} \max_{t \in I_m}{|\beta_n^{\prime\prime}(t)|} \leq  
|B_{m-1}| + \frac{\gamma_{n,2}^{-1/4}}{n} \gamma_{n,2}n^2 \leq 2n\varepsilon$$
for every $t \in I_m$ if $n$ is sufficiently large.
So
$$E_{n,\varepsilon} \subset \{t \in [0,2\pi]: |\beta_n^\prime(t)| \leq 2n\varepsilon\}$$
for every sufficiently large $n$.
Hence, by assumptions (3.2) we have
$$\text{\rm meas}(E_{n,\varepsilon}) \leq \gamma(4\varepsilon) + \gamma_n(4\varepsilon)$$
for every sufficiently large $n$. Combining this with $0 \leq {G_{n,q}(t)} \leq 1$, $t \in [0,2\pi)$, 
we obtain 
$$\left| \sum_{m}{\int_{I_m}{G_{n,q}(t) \,dt}} - K(q) \sum_{m}{\text{\rm meas}(I_m)} \right|  
\leq (\gamma(4\varepsilon) + \gamma_n(4\varepsilon)(1 + K(q))\,, \tag 3.13$$
for every sufficiently large $n$, where the summation is taken over all $m=1,2, \ldots,N$
for which $|B_{m-1}| < n\varepsilon$, and where
$\displaystyle{\lim_{\varepsilon \rightarrow 0+}{\gamma(4\varepsilon)} = 0}$ and 
$\displaystyle{\lim_{n \rightarrow \infty}{\gamma_n(4\varepsilon)} = 0}$ by assumption (3.3).
Since $\varepsilon > 0$ is arbitrary, (3.11) follows from (3.12) and (3.13). 
As we have already pointed out (3.5) follows from (3.11). 
\qed \enddemo

\demo{Proof of Lemma 3.9}
Let $\varepsilon > 0$ be fixed. Let $L_n := \gamma_{n,2}^{-1/4}$ be the same as in Lemma 3.8. 
Let the intervals $I_m$ and the numbers $A_m$ and $B_m$, $m=1,2, \ldots, N$, be the same as in the proof 
of Lemma 3.8. We define
$$F_{n,q}(t) := \left| \cos(\beta_n(t)) \right|^q\,,$$
$$\widetilde{F}_{n,q}(t) := F_{n,q}(t) \left| n^{-1} \beta_n^\prime(t) \right|^q\,,$$
$$G_{n,q}(t) := \left\{\aligned
|\cos&(A_{0} + B_{0}(t - a_{0}))|^q, \\
|\cos&(A_{1} + B_{1}(t - a_{0}))|^q, \\
& \vdots \\
|\cos&(A_{N-1} + B_{N-1}(t - a_{N-1}))|^q, \endaligned \right. \qquad
\aligned
t & \in I_1\,, \\ 
t & \in I_2\,, \\ 
& \vdots \\ 
t & \in I_N\,, \endaligned $$ 
$$H_{n,q}(t) := \left\{\aligned
|n^{-1}&B_0|^q, \\
|n^{-1}&B_1|^q, \\
& \vdots \\
|n^{-1}&B_{N-1}|^q, \endaligned \right. \qquad
\aligned
t & \in I_1\,, \\
t & \in I_2\,, \\
& \vdots \\
t & \in I_N\,, \endaligned $$
and
$$\widetilde{G}_{n,q}(t) := G_{n,q}(t)H_{n,q}(t)\,. \tag 3.14$$
Similarly to the corresponding argument in the proof of Lemma 3.8, we obtain   
$$\lim_{n \rightarrow \infty}{\sup_{t \in [0,2\pi)}{| G_{n,q}(t) - F_{n,q}(t)|}} = 0\,. \tag 3.15$$
It follows from assumption (3.4) that
$$\split \left| \left| n^{-1}\beta_n^\prime(t) \right| - \left| n^{-1}B_{m-1} \right| \right| = & 
\left| \left| n^{-1}\beta_n^\prime(t) \right| - \left| n^{-1}\beta_n^\prime(a_{m-1}) \right| \right| \cr 
\leq & \left| n^{-1}\beta_n^\prime(t) - n^{-1}\beta_n^\prime(a_{m-1}) \right| \cr 
\leq & \frac{L_n}{n} \max_{t \in I_m}{\left |n^{-1}\beta_n^{\prime\prime}(t) \right|}  
\leq \frac{\gamma_{n,2}^{-1/4}}{n} n^{-1} \gamma_{n,2}n^2 = \gamma_{n,2}^{3/4} \cr \endsplit$$
for every $t \in I_m$, where $\displaystyle{\lim_{n \rightarrow \infty}{\gamma_{n,2}^{3/4}} = 0}$. 
Hence
$$\lim_{n \rightarrow \infty}
{\sup_{t \in [0,2\pi)}{\left| H_{n,q}(t) - \left| n^{-1}\beta_n^\prime(t) \right|^q \right|}} = 0\,. \tag 3.16$$
Observe that
$$\sup_{t \in [0,2\pi)}{\left| \cos(\beta_n(t) \right|^q} \leq 1\,, \tag 3.17$$
and by the assumption (3.7) we have
$$\sup_{t \in [0,2\pi)}{\left| n^{-1} \beta_n^\prime(t) \right|^q} \leq c^q\,. \tag 3.18$$
Now (3.14), (3.16), (3.17), (3.18), and (3.14) imply
$$\lim_{n \rightarrow \infty}{\sup_{t \in [0,2\pi)}{|\widetilde{G}_{n,q}(t) - \widetilde{F}_{n,q}(t)|}} = 0\,.$$
Therefore, in order to prove (3.8), it is sufficient to prove that
$$\int_0^{2\pi}{\widetilde{G}_{n,q}(t) \,dt} \sim \frac{2\pi K(q)}{q+1}\,. \tag 3.19$$
As a special case of (3.18), we have
$$\left| n^{-1}B_{m-1} \right|^q \leq c^q\,, \qquad m =1,2, \ldots, N\,. \tag 3.20$$
If $|B_{m-1}| \geq n\varepsilon$, then (3.14), (3.20), and Lemma 3.7 give that   
$$\left| \int_{I_m}{\widetilde{G}_{n,q}(t) \,dt} - K(q)\text{\rm meas}(I_m)\left| n^{-1}B_{m-1} \right|^q \right|
\leq c^q \frac{\pi}{n\varepsilon}\,.$$
It follows from assumption (3.4) that 
$\displaystyle{\lim_{n \rightarrow \infty}{L_n} = \lim_{n \rightarrow \infty}{\gamma_{n,2}^{-1/4}} = \infty}$,
and hence
$$\split \left| \sum_{m}{\int_{I_m}{G_{n,q}(t) \,dt}} - 
K(q)\sum_{m}{\text{\rm meas}(I_m) \left| n^{-1}B_{m-1} \right|^q} \right|
& \leq N c^q \frac{\pi}{n\varepsilon} \cr
& \leq c^q \left(\frac{2\pi n}{L_n} + 1\right) \frac{\pi}{n\varepsilon} \cr 
& \leq \eta_n(\varepsilon,q)\,, \cr \endsplit \tag 3.21$$
where the summation is taken over all $m=1,2, \ldots, N$ for which $|B_{m-1}| \geq n\varepsilon$,
and where $(\eta_n(\varepsilon,q))$ is a sequence of real numbers tending to $0$  
for every fixed $\varepsilon > 0$ and $q > 0$. 
Now let
$$E_{n,\varepsilon} := \bigcup_{m: \, |B_{m-1}| \leq n\varepsilon}{I_m}\,.$$
As in the proof of Lemma 3.8 we have
$$\text{\rm meas}(E_{n,\varepsilon}) \leq 8\pi\varepsilon + \gamma_n(4\varepsilon)\,,$$
for every sufficiently large $n$. 
Combining this with (3.17) and (3.20), and recalling the definition of $\widetilde{G}_{n,q}$, we obtain
$$\left| \sum_{m}{\int_{I_m}{\widetilde{G}_{n,q}(t) \,dt}} - K(q)\sum_m{\text{\rm meas}(I_m) \left| n^{-1}B_{m-1} \right|^q} \right| 
\leq  (8\pi \varepsilon + \gamma_n(4\varepsilon))c^q(1 + K(q)) \tag 3.22$$
for every sufficiently large $n$, where the summation is taken over all $m=1,2,\ldots,N$ for which 
$|B_{m-1}| < n\varepsilon$, and where $$\lim_{n \rightarrow \infty}{\gamma_n(4\varepsilon)|} = 0$$        
by assumption (3.3). Since $\varepsilon > 0$ is arbitrary, from (3.21) and (3.22) we obtain that
$$\int_0^{2\pi}{\widetilde{G}_{n,q}(t) \,dt} \sim K(q)\int_0^{2\pi}{H_{n,q}(t) \,dt}\,. \tag 3.23$$
However, (3.16) and (3.9) imply that
$$\int_0^{2\pi}{H_{n,q}(t) \,dt} \sim n^{-q}\int_0^{2\pi}{\left| \beta_n^\prime(t) \right|^q \,dt} 
\sim \frac{2\pi}{{2}^q(q+1)}\,. \tag 3.24$$
The statement under (3.19) now follows by combining (3.23), and (3.24). 
As we have remarked before, (3.19) implies (3.8). 
\qed \enddemo

\head 4. Proofs of Theorems 2.1--2.7 \endhead  

\demo{Proof of Theorem 2.1} 
Using the notation (1.3) observe that (1.5) implies that the functions 
$$\beta_n(t) := \frac 12 (\alpha_n(t) - \alpha_n^*(t)) = \alpha_n(t) - \frac{nt}{2} - t_0 \tag 4.1$$
satisfy 
$$\beta_n^{\prime}(t) = \alpha_n^{\prime}(t) - n/2\,, \qquad t \in {\Bbb R}\,,$$
and
$$\beta_n^{\prime\prime}(t) = \alpha_n^{\prime\prime}(t)\,, \qquad t \in {\Bbb R}\,,$$
and hence Lemmas 3.1, 3.2, and 3.3 imply that the functions $\beta_n$ satisfy assumptions (3.2), (3.3), and (3.4) 
of Lemma 3.8. Hence the theorem follows from Lemmas 3.5 and 3.8. 
\qed \enddemo

\demo{Proof of Theorem 2.3}
Using the notation (1.3) let
$$P_n(e^{it}) = R_n(t) e^{i\alpha_n(t)}  \qquad  \text{and} \qquad R_n^*(t) e^{i\alpha_n^*(t)}\,,$$
where 
$$R_n(t) = |P_n(e^{it})| = |P_n^*(e^{it})| = R_n^*(t)\,,$$
and both $R_n$ and $\alpha_n$ are in $C^\infty({\Bbb R})$. 
Let $(P_n)$ be an ultraflat sequence of polynomials $P_n \in {\Cal K}_n$. We have 
$$R_n(t)^2 = n(1 + \delta_n(t))\,, \qquad t \in {\Bbb R}\,, \quad 
\lim_{n \rightarrow \infty}{\max_{t \in [0,2\pi]}{|\delta_n(t)|}} = 0\,. \tag 4.2$$ 
Let $\beta_n(t)$ be defined by (4.1). Let $\varepsilon > 0$ be fixed.
Using the cosine theorem for triangles, Lemmas 3.2, and 3.4, and (4.1), (4.2), and (1.5),   
we obtain 
$$\split & \, \left| (P_n^\prime - P_n^{*\prime}(e^{it}) \right|^2 
- \left| P_n^\prime(e^{it}) \right|^2 - \left| P_n^{*\prime})(e^{it}) \right|^2 \cr
= & \, 2((R_n(t)\alpha_n^\prime(t))^2 + (R_n^\prime(t))^2)^{1/2}
(R_n(t)(\alpha_n^{*\prime}(t))^2 + (R_n^\prime(t))^2)^{1/2} \cos(2\beta_n(t)) \cr  
= & \, 2(R_n(t)\alpha_n^\prime(t))(R_n(t)\alpha_n^{*\prime}(t))\cos(2\beta_n(t)) + \eta_n(t)n^3 \cr 
= & \, 2(R_n(t)^2(n\alpha_n^\prime(t) - \alpha_n^\prime(t)^2)\cos(2\beta_n(t)) + \eta_n(t)n^3 \cr
= & \, 2(n + \delta_n(t))(n\alpha_n^\prime(t) - \alpha_n^\prime(t)^2)\cos(2\beta_n(t)) + \eta_n(t)n^3 \cr
= & \, 2n(n\alpha_n^\prime(t) - \alpha_n^\prime(t)^2)\cos(2\beta_n(t)) + \varphi_n(t)n^3 +  \eta_n(t)n^3 \cr 
= & \, 2n\left( \frac{n^2}{4} - \beta_n^\prime(t)^2) \right) \cos(2\beta_n(t)) 
+ \varphi_n(t)n^3 +  \eta_n(t)n^3 \cr \endsplit \tag 4.3$$
with some real numbers $\eta_n(t)$ and $\varphi_n(t)$ satisfying 
$$\max_{t \in [0,2\pi]}{|\varphi_n(t) + \eta_n(t)|} < \varepsilon$$ 
for every sufficiently large $n$. Observe that
$$\int_0^{2\pi}{\left( \left| P_n^\prime(e^{it}) \right|^2 + \left| P_n^{*\prime}(e^{it}) \right|^2 \right) \, dt} 
= 4\pi \, \frac{n(n+1)(2n+1)}{6} \tag 4.4$$
by the Parseval formula. The integration by parts formula and Lemma 3.3 give that
$$\split & \, \left| \int_0^{2\pi}{\beta_n^\prime(t)^2 \cos(2\beta_n(t)) \, dt} \right| \cr 
= & \, \left| \left[ \frac 12(\sin(2\beta_n(t))\beta_n^\prime(t)) \right]_0^{2\pi} 
- \int_0^{2\pi}{\beta_n^{\prime\prime}(t)\sin(2\beta_n(t)) \, dt} \right| \cr
= & \, \left| \int_0^{2\pi}{\beta_n^{\prime\prime}(t)\sin(2\beta_n(t)) \, dt} \right| 
\leq \, \int_0^{2\pi}{|\beta_n^{\prime\prime}(t)|\, dt} = \int_0^{2\pi}{|\alpha_n^{\prime\prime}(t)|\, dt} \cr 
\leq & 2\pi \gamma_{n,2}n^2 < \varepsilon n^2 \endsplit \tag 4.5$$
for every sufficiently large $n$. Observe also that Lemma 3.8 gives that 
$$\split \int_0^{2\pi}{\cos(2\beta_n(t)) \, dt} 
= & \, \int_0^{2\pi}{\left( 2\sin^2(\beta_n(t)) - 1 \right) \, dt} 
= 2\pi(2K(2) + h_n - 1) \cr 
= & 2\pi \left(2\frac{\Gamma \left(\frac{3}{2} \right)}{\Gamma(2) \sqrt{\pi}} + h_n - 1 \right) 
= 2\pi h_n \endsplit \tag 4.6$$ 
with a sequence $(h_n)$ converging to $0$. Combining (4.3)--(4.6) we conclude 
$$\left| \frac{1}{2\pi}\int_0^{2\pi}{\left|(P_n^\prime - P_n^{*\prime})(e^{it})\right|^2 \, dt} - \frac{2n^3}{3} \right| 
\leq \varepsilon n^3$$
for every sufficiently large $n$. As $\varepsilon > 0$ is arbitrary, this finishes the proof. 
\qed \enddemo

\demo{Proof of Theorem 2.4}
Using the notation (1.3) observe that (1.5) implies that the functions $\beta_n$ defined by (4.1) satisfy
$$\beta_n^{\prime}(t) = \alpha_n^{\prime}(t) - n/2\,, \qquad t \in {\Bbb R}\,,$$
and
$$\beta_n^{\prime\prime}(t) = \alpha_n^{\prime\prime}(t)\,, \qquad t \in {\Bbb R}\,,$$
and hence Lemmas 3.1, 3.2, and 3.3 imply that the functions $\beta_n$ satisfy assumptions (3.6), (3.3), 
(3.4), and (3.7) of Lemma 3.9. Hence the theorem follows from Lemmas 3.6 and 3.9.
\qed \enddemo

\demo{Proof of Theorem 2.5}
Let $(P_n)$ be an ultraflat sequence of polynomials $P_n \in {\Cal K}_n$. Theorem 2.2 gives that 
$$\sum_{k=0}^n{\left| a_{k,n} - \overline{a}_{n-k,n} \right|^2} \sim 2n\,,$$
which is equivalent to 
$$2 \text {\rm Re} \left( \sum_{k=0}^n{a_{k,n}a_{n-k,n}} \right) = o(n)\,.$$
Now let $c \in {\Bbb C}$, $|c|=1$, and let $Q_n$ be defined by $Q_n(z) = P_n(cz)$.  Observe that $(Q_n)$ is an 
ultraflat sequence of  polynomials $Q_n \in {\Cal K}_n$ and hence
$$ \text {\rm Re} \left( \sum_{k=0}^n{c^n a_{k,n}a_{n-k,n}} \right) = o(n)\,,$$ 
and hence
$$\left| \sum_{k=0}^n{a_{k,n}a_{n-k,n}} \right| = o(n)\,.$$
\qed \enddemo

\demo{Proof of Theorem 2.6}
Let $(P_n)$ be an ultraflat sequence of polynomials $P_n \in {\Cal K}_n$. Theorem 2.3 gives that
$$\sum_{k=0}^n{k^2 \left| a_{k,n} - \overline{a}_{n-k,n} \right|^2} \sim \frac{2n^3}{3}\,,$$
which is equivalent to
$$2 \text {\rm Re} \left( \sum_{k=0}^n{k^2a_{k,n}a_{n-k,n}} \right) = o(n^3)\,.$$
Now let $c \in {\Bbb C}$, $|c|=1$, and let $Q_n$ be defined by $Q_n(z) = P_n(cz)$.  Observe that $(Q_n)$ is an
ultraflat sequence of  polynomials $Q_n \in {\Cal K}_n$ and hence
$$ \text {\rm Re} \left( \sum_{k=0}^n{c^nk^2a_{k,n}a_{n-k,n}} \right) = o(n^3)\,,$$
and hence
$$\left| \sum_{k=0}^n{k^2a_{k,n}a_{n-k,n}} \right| = o(n^3)\,.$$
\qed \enddemo

\demo{Proof of Theorem 2.7.}
Using Lemma 3.5* and then applying Lemma 3.8 with $\beta_{2n}$ defined by $\beta_{2n}(t) := \alpha_n(t) + \pi/2$ 
we obtain the first asymptotic formula of the theorem. Using Lemma 3.6* and then applying Lemma 3.9 with $\beta_{2n}$ 
defined by $\beta_{2n}(t) := \alpha_n(t) + \pi/2$ we obtain the second asymptotic formula of the theorem.
\qed \enddemo

\enddocument